\documentclass[sigconf]{acmart}

\usepackage{booktabs} 

\setcopyright{rightsretained}

\usepackage[utf8]{inputenc}

\usepackage{amssymb,amsmath,url,rotating,ifthen,algorithm,algorithmicx,epsfig,multirow,array,color,multicol,graphicx,float,hhline,calc,algpseudocode}
\usepackage{todonotes} 
\usepackage{bold-extra,ulem}
\usepackage{enumitem}
\setlistdepth{9}
\setlist[itemize,1]{label=$\bullet$}
\setlist[itemize,2]{label=$-$}
\setlist[itemize,3]{label=$\ast$}
\setlist[itemize,4]{label=$\cdot$}
\setlist[itemize,5]{label=$\diamond$}
\renewlist{itemize}{itemize}{5}

\newcommand{\ot}[1]{\textcolor{black}{#1}}
\newcommand{\otc}[1]{\textcolor{black}{#1}}

\algdef{SE}[PROCEDURE]{Procedure}{EndProcedure}%
[2]{\algorithmicprocedure\ \textproc{#1}\ifthenelse{\equal{#2}{}}{}{(#2)}}%
{\algorithmicend\ \algorithmicprocedure}%
   
\algdef{SE}[FUNCTION]{Function}{EndFunction}%
[2]{\algorithmicfunction\ \textproc{#1}\ifthenelse{\equal{#2}{}}{}{(#2)}}%
{\algorithmicend\ \algorithmicfunction}%

\usepackage{hyperref}
\hypersetup{
    colorlinks,
    citecolor=black,
    filecolor=black,
    linkcolor=blue,
    urlcolor=black}

\def\CQFD{\fbox{}}
\def\qed{\fbox{}}

\def\E{{\mathbb E}}

\def\N{\mathcal{N}}   

\def\e{\epsilon }

\def\R{{\mathbb R}}

\begin{document}
\title{Population Control meets Doob's Martingale Theorems: the Noise-free Multimodal Case}
\author{
Marie-Liesse Cauwet,
Olivier Teytaud
}

\begin{abstract}
We study a test-based population size adaptation (TBPSA) method, inspired from population control, in the noise-free multimodal case.
In the noisy setting, TBPSA usually recommends, at the end of the run, the center of the Gaussian as an approximation of the optimum.
We show that combined with a more naive recommendation, namely recommending the visited point which had the best fitness value so far, TBPSA is also powerful in the noise-free multimodal context.
	We demonstrate this experimentally and explore this mechanism theoretically: we prove that TBPSA is able to escape plateaus with probability one in spite of the fact that it can converge to local minima. This leads to an algorithm effective in the multimodal setting without \otc{resorting} to a random restart from scratch.


\def\reremoveme{Deadlines:
\begin{itemize}
    \item  16th International Conference on Parallel Problem Solving from Nature (PPSN XVI)	05/09/2020	01/04/2020
\item Genetic and Evolutionary Computation Conference (GECCO 2020)	08/07/2020		30/01/2020
\end{itemize}
}
\end{abstract}

\maketitle


\section{Introduction}
\subsection{Population control}
Population control has been
proposed in \cite{beyernoise}  and adapted in \cite{nevergrad} under the name TBPSA (slightly different from the original population control) with great successes in noisy optimization. Consistently with \cite{sandranoiseuhcma}, TBPSA breaks the barrier of a simple regret 
\begin{equation}
simple\ regret = O(1/\sqrt{number\ of\ evaluations}).\label{slow}
\end{equation}by doing steps in the recommendation smaller than the steps in the exploration; this is the so-called ``mutate large, inherit small'' paradigm~\cite{mlis}, i.e. we must explore ``far'' from the approximate optimum for reaching simple regret 
\begin{equation}
simple\ regret = O(1/number\ of\ evaluations).\label{fast}
\end{equation}
\subsection{Multimodal optimization}
In the present paper we 
consider the application of the same approach to multimodal noise-free optimization.
Multimodal optimization can be the search for a global optimum, in a context made difficult by the presence of many local optima.
In other contexts, it refers to the search for diverse global optima.
We show, in the case of a plateau, that though population control does not necessarily lead to large step-sizes it will nonetheless sample far thanks to a preserved diversity (i.e. no fast decrease of the step-size to zero) so that the increasing population will mechanically provide points far enough for escaping a local minimum.

\subsection{\otc{Known} Convergence Results}
\subsubsection{Convergence rate in the noise-free case.}
It is known~\cite{Rechenberg,beyerbook,youheifoga2019} that evolutionary algorithms converge in many cases to the optimum $x^*$ with rate $\log\|x_n-x^*\|  =-\Omega(n)$. 
\begin{figure*}[t]
\centering
\includegraphics[width=.3\textwidth]{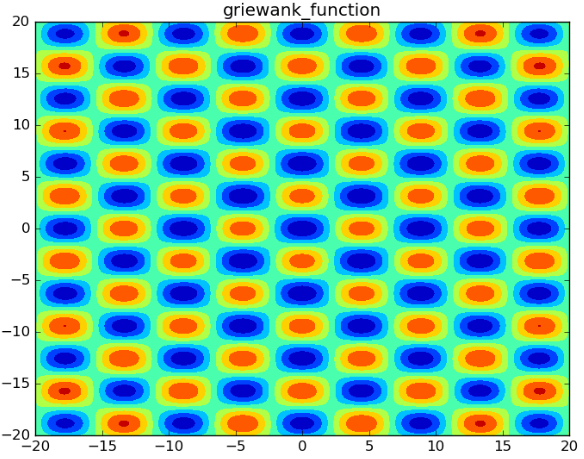}
\includegraphics[width=.3\textwidth]{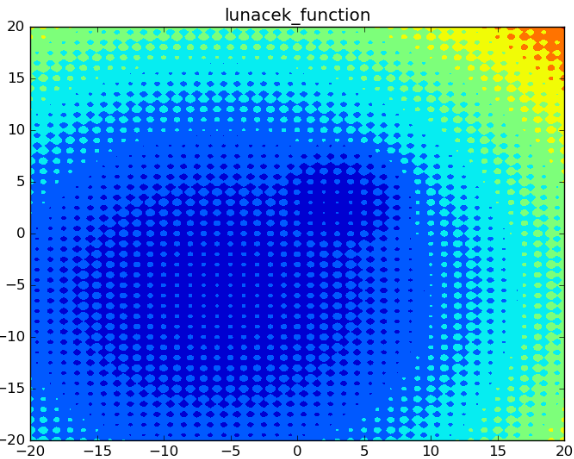}
\includegraphics[width=.3\textwidth]{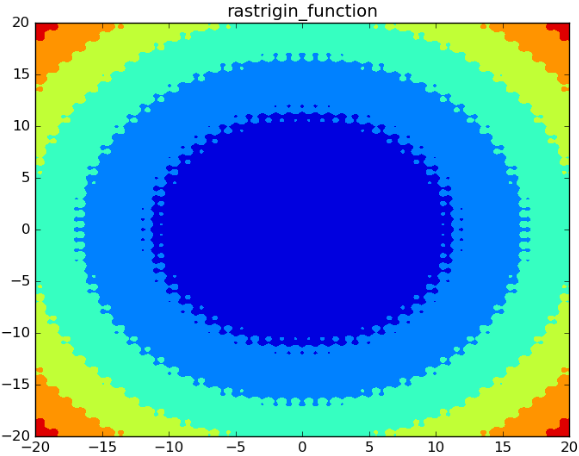}
\includegraphics[width=.3\textwidth]{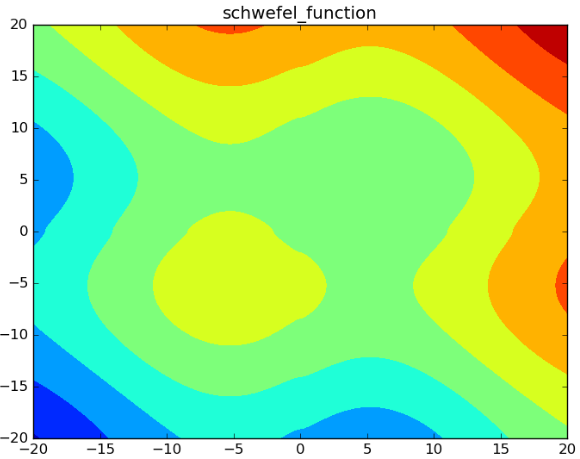}
\includegraphics[width=.3\textwidth]{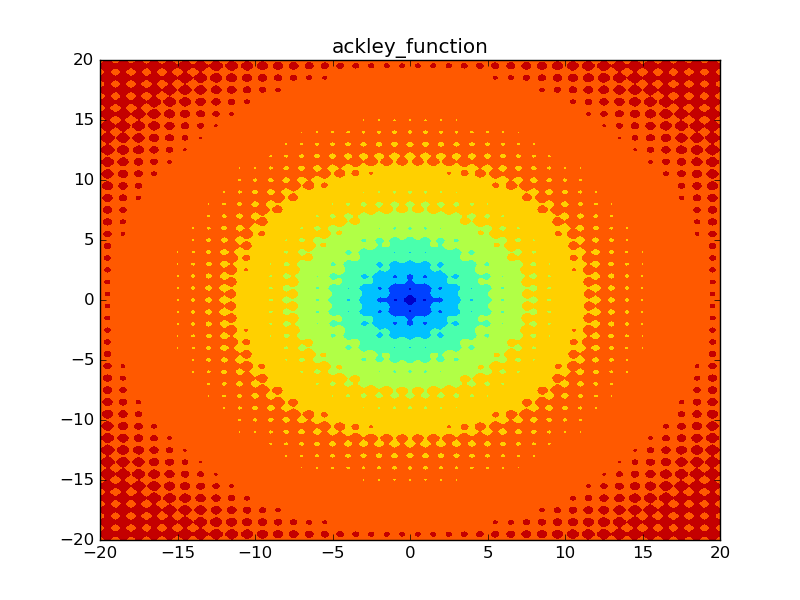}
\caption{\label{plots}Level sets of some of our test functions. The optimum is in the bluest parts.}
\end{figure*}
\begin{figure*}[t]
    \centering
    \includegraphics[trim={0 10 0 90}, clip, width=.8\textwidth]{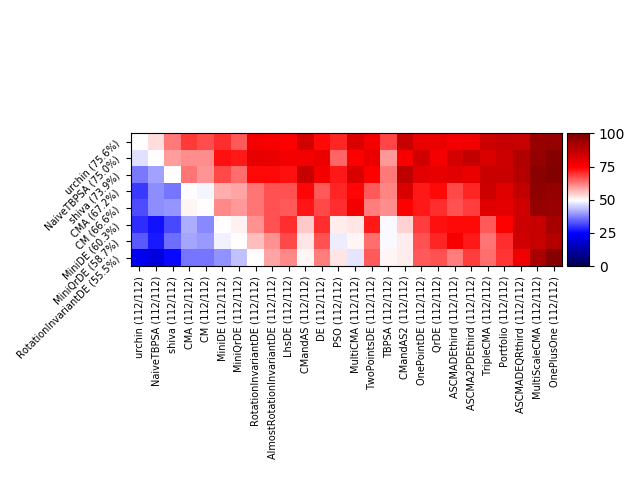}
     \includegraphics[trim={0 10 0 90}, clip,width=.8\textwidth]{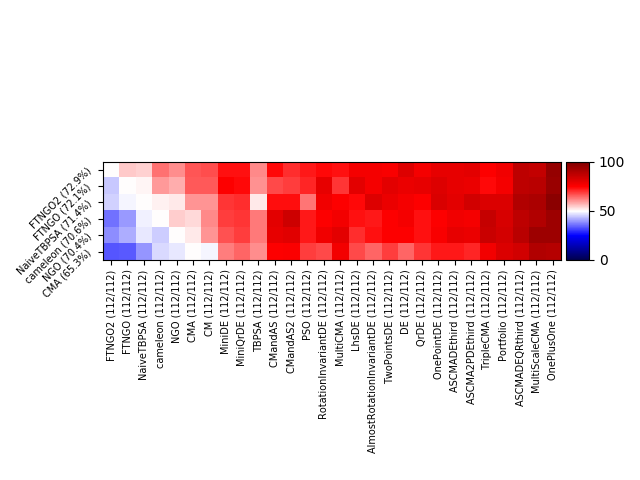}
    \caption{Parallel  multimodal experiment in Nevergrad (termed ``paramultimodal'' in \cite{nevergrad}) at two distinct times: the list of algorithms varies, NaiveTBPSA performs well in both. Methods are ranked by score (best at the top row, best on the left hand side column). Only the best rows are presented. The score is the average frequency at which method $x$ outperformed other methods over instances of test functions. All other methods in the top (JNGO, FTNGO, cameleon, octopus, NGO, Shiva, Urchin) are combinations of algorithms based on algorithm selection: \textit{they incorporate our own method naiveTBPSA} - all other methods are significantly weaker, CMA ranking best among methods not using NaiveTBPSA. This experiment considers 1000 concurrent function evaluations (parallelism), budget in $\{3000, 10000, 30000, 100000\}$. We refer to \cite{nevergrad} for a detailed description of all algorithms involved in the comparison and for the detailed setup. The objective functions are "Hm", "Rastrigin", "Griewank", "Rosenbrock", "Ackley", "Lunacek", "DeceptiveMultimodal", which are all either well known\ot{, or recent (namely Hm and DeceptiveMultimodal) but open sourced in \cite{nevergrad}.}}
    \label{paramulti}
\end{figure*}
\begin{figure*}[t]
    \centering
    \includegraphics[width=.7\textwidth]{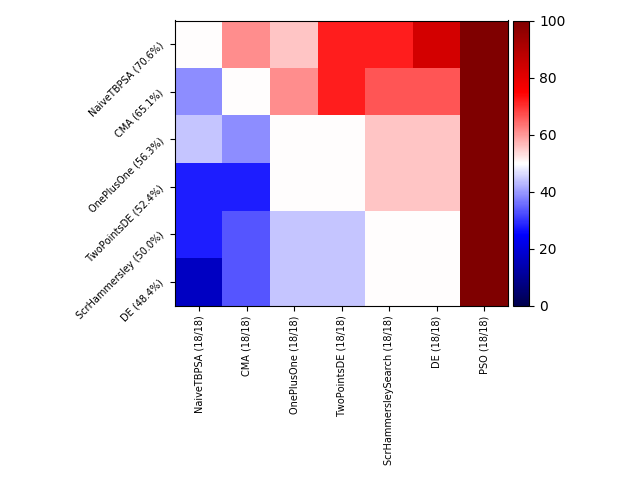}
    \includegraphics[width=.7\textwidth]{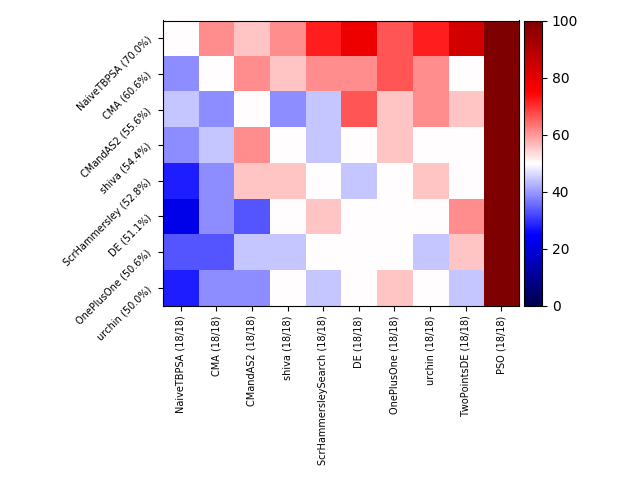}
    \caption{Parallel experiment in Nevergrad at two distinct dates with various algorithms present in that experiment at these distinct moments: this was not our goal, but we see that our code performs quite well in the parallel setting, in this experiment termed ``parallel'' (number of parallel evaluations equal to 20\% of the budget, budget in $\{30,100,3000\}$, objective function in $\{ Sphere, Rastrigin, Cigar\}$) and defined in \cite{nevergrad}.}
    \label{para}
\end{figure*}
\begin{figure*}[t]
    \centering
    \includegraphics[trim={0 10 0 90}, clip,width=.8\textwidth]{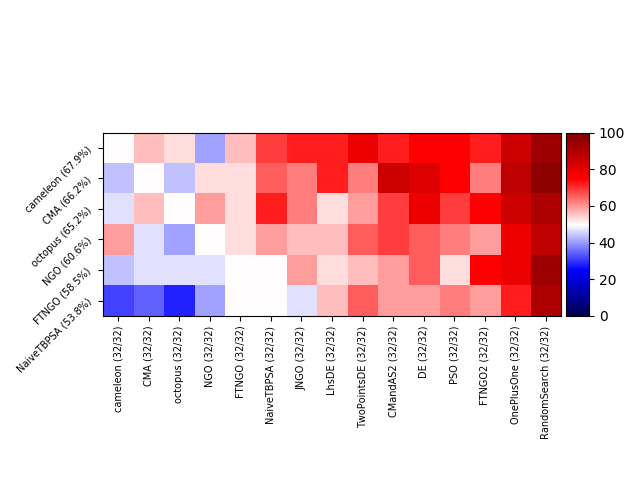}
    \includegraphics[trim={0 10 0 90}, clip,width=.8\textwidth]{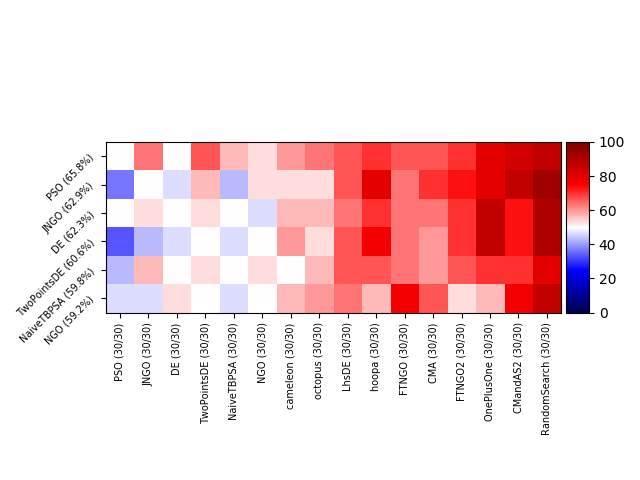}    \caption{Multiobjective (top, set of problems with 2 or 3 objectives, budget from $100$ to $2900$, sequential or $100$ parallel workers) and Manyobjective (bottom, 6 objectives, budget from $100$ to $5900$, sequential or $100$ parallel workers) experiment in Nevergrad (code available in \cite{nevergrad}). This experiment uses the hypervolume indicator, i.e. the obtained method problem is multimodal. NaiveTBPSA performs better than PSO and all DE variants in the multiobjective case, but worse than CMA. It performs better than CMA and some variants of DE in the manyobjective case, but worse than some variants of DE and PSO. NaiveTBPSA has the best average rank over the two methods among simple methods. Algorithms such as Octopus, Cameleron, and names including ``NGO'' are combinations of algorithms designed independently of the present work: they are actually built on top of NaiveTBPSA and many others combined through an algorithm selection method.}
    \label{hyperv}
\end{figure*}
\subsubsection{Convergence in the multimodal case.}
\otc{Multimodal optimization can be tackled in a number of ways~\cite{multimodal}.
A wide research field in multimodal optimization consists in niching methods\cite{holland37,niching2}, clearing\cite{clearing}, sharing\cite{sharing,sharing2}. As opposed to restart algorithms, these methods preserve diversity {\em{during}} the optimization run. At the cost of slower local convergence because of several local convergence runs simultaneously, these methods are more valuable when we want to find several optima (as opposed to just finding one of the global optima) or when we have to be parallel. For finding a single global optimum, restarts remain the method of choice and dominates e.g. \cite{HAFR2012RPBBOBES}.
When an algorithm stagnates, restart methods assume that it is stuck in a local optimum, and launches another run from a random initial point\cite{auger2005restart}. \cite{qrrs} shows that quasi-random restarts are faster than random restarts. An improvement consists in adding a bandit for choosing between independent runs\cite{banditqrrs}.
Theoretical investigations of restarts do exist.  \cite{ast} shows that random diversification can be combined with evolution strategies for having both global convergence and reasonably fast rates. More precisely, it shows that a local linear convergence result, or a convergence faster than linear, is preserved when we use a random diversification: this combines the best of both worlds, almost sure convergence on a wide range of functions as in random search, and fast local convergence. This solution for multimodal convergence, however, is computationally expensive non-asymptotically as the multimodalities are handled through random search.}

\subsection{Optimization Algorithms}

\def\removed{
For presenting optimization algorithms, we distinguish exploration and recommendation. Exploration means the way we choose which point is going to be evaluated next; whereas recommendation means the way we choose the approximation of the optimum returned at the end of the optimization run.
Unless stated otherwise, algorithms use for recommendation the point with best obtained fitness.
A more detailed list of algorithm is as follows:
\begin{itemize}
\item Random search: exploration by randomly drawing a Gaussian point.
\item Quasirandom point and design of experiments: we use Latin Hypercube Sampling\cite{eglajs,mckay}, Halton\cite{halton}, Hammersley\cite{hammersley}, and scrambled versions thereof\cite{atanasov,permut1,permut2,permut3,owen}. These algorithms are usually presented in terms of sampling the hypercube $[0,1]^d$; we convert to unbounded domains by applying the inverse cumulative distribution function of the Gaussian for each coordinate.
\item toto
\end{itemize}
}
We present below TBPSA, the Test-Based Population-Size Adaptation method from \cite{nevergrad}. We will then present our counterpart, \ot{termed NaiveTBPSA,} equipped with a different recommendation method.
We also use many algorithms from \cite{nevergrad} in our experiments. Readers unfamiliar with evolution strategies (ES) are referred to e.g. \cite{beyerbook}.
\subsubsection{TBPSA: population control for noise management.}
We use a TBPSA self-adaptive $(\mu/\mu,\lambda)$-ES, implemented in Nevergrad, and strongly inspired from \cite{beyernoise}. More precisely, each point is a pair $(x_i,\sigma_i)$ where $x_i$ is a candidate solution and the step-size $\sigma_i$ is a positive number. At each generation we compare their fitness values, select the $\mu$ best, average the selected candidates for choosing the next parent $x$, log-average their step-sizes for choosing the next parent step-size $\sigma$, and generate $\lambda$ points by $\sigma_i=\sigma \exp(\N)$ and $x_i=x+\sigma_i \N_d$ for $1\leq i \leq \lambda$. Consistently with \cite{beyernoise}, we multiply both $\lambda$ and $\mu$ by 2 if the $\lambda$ most recent and the $\lambda$ oldest out of the $5\lambda$ last points have no statistical difference. \otc{No statistical differnce means that} the differences between averages is not 2 standard deviations apart - we multiply $\lambda$ and $\mu$ by $2^{-\frac14}$ otherwise, without ever decreasing below the initial value or below the degree of parallelism requested by the application. Initialization: $\lambda=4d$, $\mu=d$, $\sigma=1/\sqrt{d}$. In classical applications of population control, this algorithm uses the last parent as a recommmendation. These two features (test-based population size, recommendation equal to the parent) are critical for reaching a convergence rate ``simple regret $= O(1/\mbox{budget})$'' rather than $1/\sqrt{\mbox{budget}}$ in the noisy case. 
\subsubsection{Naive TBPSA: TBPSA with best so far recommendation. }

We define  naiveTBPSA-$(\mu/\mu,\lambda)$-ES: the same as above, but using the same recommendation method as most algorithms in the noise-free context, namely the best visited point from the point of view of their fitness when visited.


\def\totoremoved{\begin{itemize}\item We use differential evolution\cite{de}, in several flavors:
\begin{itemize}
   \item Default algorithm, with population size $5\times $, Curr-To-Best, $CR=.8$, $F=\frac12$.
   \item Berthier's variant\cite{berthierDE}, with Curr-To-Best, $F1=.8$, $F2=.8$.\todo{which name in the xps, and do we keep this}
   \item OnePointDE: the same with 1 point crossover as in \cite{1ptcrossover}.
   \item TwoPointsDE: the same but with 2 points crossover as in \cite{2ptcrossover}.
\end{itemize}
\item \todo{do we keep this ? NAS, a refactor version of the optimizer presented in \cite{nas}. TODO differences}
TODO this was optimized specifically for tf.hub or ... ?
TODO running the original algorithm NAS as well.
\item Random Search Plus Middle Point just takes care to have at least one point in the very middle as a first exploration and then applies classical random search.
\item The one-plus-one algorithm comes from \cite{rechenberg73} and uses step-size multiplied by $2$ in case of success and by $2^{-1/4}$ otherwise. The step-size is initialized at $1/\sqrt{dimension}$.
\todo{Do we keep this ? We also use Vizier's random search, converted to a domain $\R^d$ by applying a linear mapping to $[-1,1]^d$ and then the arctanh transformation for each coordinate. TODO do we have this, do we keep it ?}
\end{itemize}}

\section{Can we escape local minima ?}
A key question for an optimization algorithm applied to multimodal objective functions is its ability to escape local minima. We distinguish two cases:
\begin{itemize}
\item Convex local optima for which the evolution strategy converges log-linearly. We show that adding an exponential increase of the population size is not enough for escaping such local optima.
\item Plateaus. Then we show that the algorithm escapes plateaus, by stabilizing the step-size and increasing the population size.
\end{itemize}

\subsection{Test-based population increase can not escape local minima if evolution strategies converge fast in this local optimum}

\begin{definition}
An algorithm $A$ is Gaussian-evolutionary if, with $\tau>0$, $\N_{n,i}$ a $d$-dimensional standard normal random variables and $\N'_{n,i}$ a 1-dimensional standard normal random variable (all independent):
\begin{eqnarray}
\forall n,\lambda_n&\leq& M\lambda_{n-1}\label{e1}\\
\forall n, i\leq \lambda_n, \sigma_{n,i}&=&\sigma_n\exp(\tau \N'_{n,i} )\label{e1b}\\
\forall n, i\leq \lambda_n, x_{n,i}&=&x_n+\sigma_{n,i} \N_{n,i}\label{e2}\\
\forall n, i\leq\lambda_n, x_{n+1},\sigma_{n,i}\mbox{ and } \sigma_{n+1}& & \mbox{ depend only on the }\nonumber\\
\mbox{ $(x_n, \sigma_n, (x_{n,i}, \sigma_{n,i},f_{x_{n,i}})_{i\leq \lambda_n})$  }& &\label{e3}
\end{eqnarray}
\end{definition}
$\tau=0$ is possible and leads to a single step-size i.e. $\forall n,i, \sigma_{n,i}=\sigma_n$.

\begin{definition}
We define a property $H1$ as follows:
$$H1(K): \forall n, ||x_n|| \leq K\mbox{ and }\sigma_{n,i}\leq K \exp(-n/K).$$
\end{definition}

Property $H1$ means that the step size decreases exponentially and that points stay in a bounded set. Let us see why it is actually quite usual that $H1$ holds with some probability for some values of $K$.
For the first part, namely the existence of $K$ such that $\forall n, ||x_n||\leq K$, it is straightforward that it holds in the following case:
\begin{itemize}
    \item an elitist strategy;
    \item a coercive objective function, i.e. $\forall K>0, \exists B, \forall x, ||x||>B\Rightarrow f(x)>K$.
\end{itemize} 
And for the second part, namely the exponential decrease of the step-size, \cite{upu} and \cite{youheifoga2019} show such an exponential convergence for wide ranges of functions.
    
Let us define a property of local convergence. This property will be used as an assumption in our results.

\begin{definition}
An algorithm $A$ is locally convergent on a function $f$, denoted $LC(A,f)$ if
\begin{equation}\lim_{K\to\infty} P(H1(K))=1.\label{ha}
\end{equation}
\end{definition}

Remark: this is actually not a definition of convergence towards an optimum. This just means that we stay in a bounded neighborhood, and that 
the step-size decreases quickly. This is in fact a consequence of convergence, not a convergence. We just use this definition because it is weaker than a classical convergence to the optimum, and enough for our purpose. By using a weaker assumption, we strenghten our result; our theorem holds for this weak notion of local convergence, so a fortiori it holds for any stricter notion of local convergence.

\begin{definition}
A Gaussian evolutionary algorithm approaches the optimum on an objective function $f$ if, almost surely, for all $\e>0$, there is $n$ and $i\leq \lambda_n$ such that $\exists x^* \in \arg\min f, ||x_n-x^*||\leq \e \mbox{ or }||x_{n,i}-x^*||\leq \e$.
\end{definition}

Remarks:
\begin{itemize}
\item This definition could be extended to non-evolutionary or non-Gaussian algorithms.
\item If an algorithm does not approach the optimum, then the hitting time, for sufficiently small precision, is infinite.
\end{itemize}
\begin{theorem}[exponentially increasing the population size is not enough for escaping local minima]\label{thmlm}
Consider an algorithm such that Eqs \ref{e1}-\ref{e3} hold. We assume that 
for the fitness function $f$, $LC(A,f)$ holds.
Then, for any $\e>0$, there exists $K'=K'(\e)>0$ such that with
probability at least $1-\e$, $\forall n,i,\ ||x_{n,i}|| \leq  K'$.
\end{theorem}


This theorem can be \otc{rephrased }as follows, for showing that it implies that we can not escape local minima. 
\begin{corollary}[Corollary of Theorem \ref{thmlm}: a locally convergent Gaussian evolutionary algorithm does not escape local minima]\label{maincr}
Consider $M\in\R$, and an algorithm $A$ such that Eqs \ref{e1}-\ref{e3} hold, and
$LC(A,f)$.

Consider $\e>0$. Let $K'(\e)$ be as in Theorem \ref{thmlm}.
Consider an objective function $g$ such that $\forall x \in B(0,K'), g(x)=f(x)$ and $\arg\min g\not\in B(0,K')$. Then, with probability at least $1-\e$,
none of the $x_{n,i}$ or the $x_n$ is outside $B(0,K')$ and therefore the algorithm does not approach the optimum.
\end{corollary}



{\bf{Proof of the theorem:}}
Using Eq. \ref{ha}, let us choose $K$ such that 
\begin{equation}
P(H1(K))>1-\e/2.\label{lm1}
\end{equation}

Let us define $H2(K,K')$, for $K'>K$:
$$H2(K,K'): \forall n>0, 1\leq i\leq \lambda_n, ||{ \N}_{n,i}||\leq \frac{K'-K}{K}\exp(n/K)$$
$$\mbox{ where the }\N_{n,i}$${ are independent standard d-dimensional Gaussian random variables.}

Algebra yields:
$\forall K,K'$, 
\begin{eqnarray}
\mbox{if $H1(K)$ and $H2(K,K')$,}\nonumber\\
\mbox{then for all $n,i$ we have $||x_{n,i}||\leq K'$.}\label{lm2}
\end{eqnarray}

Using the bound $P(||\N_{n,i}|| <t)\geq 1-\alpha \exp(-\beta t)$ for some $\alpha>0$ and $\beta>0$ and $t_n=\frac{K'-K}K\exp(\frac nK)$, we get $P(\sup_{i\leq \lambda_n} ||\N_{n,i}|| <t_n)$

\begin{eqnarray*}
&\geq& 1-\alpha\lambda_0M^n \exp(-\beta t_n)\mbox{ using Eq. \ref{e1}}\\
   & \geq & 1-\alpha\lambda_0 \exp\left(n\log(M) - \beta \frac{K'-K}{K} \exp(\frac nK)\right)\end{eqnarray*}
i.e.
$P(\forall n,\sup_{i\leq \lambda_n} ||\N_{n,i}|| <t_n)$  $$\geq  1-\alpha\lambda_0\sum_n \exp\left(n\log(M) - \beta \frac{K'-K}{K} \exp(\frac nK)\right).$$

For $K'$ large enough, the right hand side is arbitrarily close to $1$.
So we get:
\begin{equation}
\forall \e, \exists h:\R\to\R, P(H2(K,K'))\geq 1-\e/2\mbox{ if }K'>h(K).\label{lm3}
\end{equation}

Then, for such $K$ and $K'$, the probability that none of the $x_{n,i}$ verifies $||x_{n,i}||>K'$ is
\begin{eqnarray*}
&\geq & P(H2(K,K')\mbox{ and }H1(K))\mbox{  by Eq.\ref{lm2}}\\
&\geq & 1-\e/2 + P(H1(K))-1 \mbox{ by Eq. \ref{lm3}}\\
&\geq & 1-\e/2 + (1-\e/2)-1\mbox{ by Eq. \ref{lm1}}\\
&\geq & 1-\e
\end{eqnarray*}
hence the expected result.
\CQFD

\subsection{Test-based population increase can escape plateaus}

There are several solutions for escaping plateaus:
\begin{itemize}
\item increasing the step-size;
\item maintaining the diversity, i.e. ensuring that the step-size does not decrease to zero;
\item adding random diversification over the domain as in \cite{ast}.
\end{itemize}
We show below that TBPSA successfully escapes plateaus by maintaining the diversity.

We consider the following algorithm, directly inspired (though not completely equal) from \cite{beyernoise}:
\begin{eqnarray}
s_{n,i}&=& \exp(\N)\\
\sigma_{n,i}&=& \sigma_n s_{n,i}\label{eq1}\\
z_{n,i}&=&\N_d\\
x_{n,i}&=&x_n + \sigma_{n,i}z_{n,i}\\
y_{n,i}&=&fitness(x_{n,i})\\
I=I_n&=& \mbox{indices of the $\mu$ best $y_{n,i}$}\\
\sigma_{n+1}&=&\exp(\mbox{average of the $\log \sigma_{n,i}$ for $i\in I$})\label{eq6}
\end{eqnarray}
where $i$, unless stated otherwise, ranges over $\{1,\dots,\lambda_n\}$ and $\N_d$ denotes an independent $d$-dimensional Gaussian standard random variable and $\N$ an independent Gaussian random variable.

\ot{We did not specify how $x_{n+1}$ is defined; our result is independent of this.}

We assume randomly broken ties.

\ot{We assume $\mu/\lambda$ lower-bounded.}

For convenience, let us note $z'_{n,i}=\sigma_{n,i}z_{n,i}/\sigma_n$;
$z'_{n,i}$ is distributed as $\exp(\N)\N_d$ and
\begin{equation}
x_{n,i}=x_n+\sigma_{n} z'_{n,i}.\label{stuce}
\end{equation}

\begin{lemma}\label{martingale}
In case of random selection,
$\E \log \sigma_{n+1}=\E\log \sigma_n$.
\end{lemma}

{\bf{Proof:}} By Eqs. \ref{eq1} and \ref{eq6}.\qed

\begin{lemma}\label{finitevar}
Let us assume that $\lambda$ is multiplied by $2$ at each iteration
and let us assume random selection. Then
the supremum over $n$ of the variance of $\log \sigma_{n+1}$ is finite.
\end{lemma}

{\bf{Proof:}} Eq \ref{eq6} is an average over $\lambda$ independent points, $\lambda$ increases exponentially with the iteration index $n$, so the variance of $\log \sigma_{n+1}$ is the variance of $\log \sigma_n$ plus $Var(\N)/2^n$ hence $Var(\log(\sigma_n))+1/2^n$. The total variance is therefore bounded by $2$.\qed

\begin{lemma}\label{doob}
 With probability $1$, $\log \sigma_n$ converges to a finite limit $\log \sigma^*$.
\end{lemma}

{\bf{Proof:}} By Doob's first martingale theorem~\cite{wikidoob}, using Lemmas \ref{martingale} and \ref{finitevar}.\qed

\begin{lemma}\label{key}
With probability $1$, $\sup_{n\leq N} \sup_{i\leq \lambda_n} ||x_{n,i}-x_n||$ goes to infinity as $N\to\infty$.
\end{lemma}

{\bf{Proof:}}
With probability $1$,
 $\sigma_n$ converges to a finite limit (by Lemma \ref{doob}) and the supremum of the $z'_{n,i}$ for $n\leq N,i\leq \lambda_n$ converges to infinity as $N\to\infty$; so Eq. \ref{stuce} concludes.
\qed

\begin{theorem}
\ot{
Define $S$ a subset of the domain on which $f$ is constant.
Then with probability $1$,  there exists $n>0,i\leq \lambda_n$ such that either $x_n\not \in S$ or  $x_{n,i}\not \in S$.}
\end{theorem}


{\bf{Proof:}}
Consider $E_N$ the event ``for all $n<N$ and $i$, $x_n$ and $x_{n,i}$ are all in $S$''.
We have two random components:
\begin{itemize}
    \item the $z_{n,i}$ and $s_{n,i}$;
    \item a random ranking, used for breaking ties.
\end{itemize}
As often done in theoretical analyses, we consider what would happen if the ranking was random instead of being based on the objective function. This corresponds to the case of random selection, i.e. plateaus and randomly broken ties.

Let us consider $x'_n$ and $x'_{n,i}$, counterpart, for the random selection case, of $x_n$ and $x_{n,i}$. Then, $x_n, x'_n, x_{n,i}, x'_{n,i}$ live in the same domain and are defined for the same universe.
We then consider $E_{rs,N}$, which is $E_N$ under random selection (using $x'$ instead of $x$).
And we consider $E_{f,N}$, which is $E_N$ under selection with $f$.
First, consider $x_n$ and $x_{n,i}$ under random selection. $f$ has no impact on what is going on.
Then, Lemma \ref{key} applies.
Therefore, with probability $1$, there is $n_0$ and $i_0$ such that $||x'_{n_0,i_0}-x'_{n_0}|| >  diameter(S)$. This implies that 
\begin{equation}
    P(\cap_{N\geq 0} E_{rs,N})=0.\label{zebeauzero}
\end{equation}
If $E_{f,N}$ holds, then for all $n\neq N$ and $i\leq \lambda_n$, the probability distribution of $(x_n,(x_{n,i})_{i\leq \lambda_n})$
 and $(x'_n,(x'_{n,i})_{i\leq \lambda_n})$ coincide, if $E_{f,n}$ holds.
This means that $E_{f,N}$ implies $E_{rs,N}$.
In particular $P(E_{f,N})\leq P(E_{rs,N})$. This and Eq. \ref{zebeauzero}
imply that $P(\cap_N E_{f,N})=0$.
\qed






\def\rrrr{
\section{Practical recommendations}
In short, we have shown that TBPSA
will not escape local minima unless they are plateaus. With plateaus, the population size increases exponentially and the step-size converges to a finite value, leading to escaping the local minimum at some point.

For local minima, we might trust the limited machine precision for minima to become plateaus. More precisely,  $f(x)$ converges to $f(x^*)$ quadratically faster than $x$ converges to $x^*$, hence we reach machine precision over fitness values earlier than in the domain.
If we want to think of machines with arbitrary precision, we can consider $x\mapsto \lfloor f(x)/\epsilon\rfloor\times\epsilon$.\todo{does this exist in the literature ?}

}

\section{Experimental results}
Fig. \ref{paramulti} and \ref{para} present results in the \otc{parallel multimodal} and multimodal setting respectively.
Fig. \ref{hyperv} presents results in the multiobjective and manyobjective setting of Nevergrad. These test cases correspond to the use of the hypervolume indicator for converting the multiobjective setting into the monoobjective case.
Figure \ref{plots} presents some of our objective functions. The detailed experimental setup is the default one in \cite{nevergrad}, all the code is public.

\def\zebigmessouille{
\subsection{Algorithms}
For presenting optimization algorithms, we must distinguish exploration and recommendation. {\em{Exploration}} means the way we choose which point is going to be evaluated next (i.e. which architecture is going to be trained); whereas {\em{recommendation}} means the way we choose the approximation of the optimum returned at the end of the optimization run. Unless stated otherwise (as for the noise-managing evolution strategies), algorithms use for recommendation the point with best obtained fitness.
The {\em{budget}} is the number of exploration points before a recommendation is issued.

We distinguish the following families of optimization algorithms, based on the way they encode their current knowledge of the objective function:
\begin{itemize}
\item Bayesian optimization \& mathematical programming methods. This includes Bayesian Optimization\cite{ego}, including entropy-based parallelization\cite{entropybo}; these methods typically have a probabilistic model of objective functions and condition it to the observed fitness values. Diversity is ensured by entropy methods.
We use Vizier's default tool for optimization as presented in \cite{vizier}; this is based on Gaussian processes (Bayesian optimization) with entropy. TODO more on Vizier default
\item Evolution strategies (ES). Evolution strategies include  methods updating a Gaussian distribution at each iteration, such as CMA-ES\cite{vizier}, CMSA\cite{CMSA}, and Evolution Strategies in general\cite{beyerbook} including (great) ideas from the recent \cite{beyernoise}. 
	We use two self-adaptive $(\mu/\mu,\lambda)$-ES strongly inspired from \cite{beyernoise}. More precisely, each point is a pair $(x_i,\sigma_i)$ where $x_i$ is a vectorized architecture and the step-size $\sigma_i$ is a positive number. At each generation we compare their fitness values (validation errors), select the $\mu$ best, average the selected vectorized architectures for choosing the next parent $x$, log-average their step-sizes for choosing the next parent step-size $\sigma$; and generate $\lambda$ points by $\sigma_i=\sigma \exp(\N)$ and $x_i=x+\sigma \N_d$. Consistently with \cite{beyernoise}, \otc{TBPSA multiplies} both $\lambda$ and $\mu$ by 2 if the $\lambda$ most recent and the $\lambda$ oldest out of the $5\lambda$ last points have no statistical difference in the sense that the differences between averages is not 2 standard deviations apart - we multiply $\lambda$ and $\mu$ by $2^{-\frac14}$ otherwise, without ever decreasing below the initial value. Initialization: $\lambda=4d$, $\mu=d$, $\sigma=1/\sqrt{d}$. In the noise-managing-$(\mu/\mu,\lambda)$-ES, this algorithm uses the last parent as a recommendation. These two features (test-based population size, recommendation equal to the parent) are critical for reaching a convergence rate ``simple regret $= O(1/\mbox{budget})$'' rather than $1/\sqrt{\mbox{budget}}$. 
\item The one-plus-one ES algorithm comes from \cite{onefifth} and uses step-size multiplied by $2$ in case of success and by $2^{-1/4}$ otherwise. The step-size is initialized at $1/\sqrt{dimension}$. 
\item Other methods modeling the current knowledge of the objective function by a set of good candidates and mixing them by crossovers (and not only averaging, contrarily to evolution strategies). This includes Differential Evolution (DE\cite{de}). 

\item We use differential evolution\cite{de}, in several flavors:
   the default algorithm, with population size $5$ times the dimension, Rand/1, $CR=.8$, $F=\frac12$; we also use a tuned variant (TDE for Tuned-DE) \cite{debt}, with Curr-To-Best, $F1=.8$, $F2=.8$; and two variants of this TDE method, combined with 1 point crossover and 2 points crossover\cite{holland37}. These versions are supposed to be somehow modular; and the order of variables has an impact when we use k-points crossover.
TDE with adapted crossover rate $CR=1$ corresponds to a rotationally invariant DE; we refer to it as TRotInvDE.

\item Design of experiments, also known as one-shot optimization methods; this includes random search, quasi-random search, latin hypercube sampling. These methods are weaker in terms of performance for a fixed number of evaluations; however, they can greatly benefit from pruning, be quite robust, and are highly parallel. Random search: exploration by randomly drawing a Gaussian point. Random Search Plus Middle Point is similar but takes care to have the first point in the very middle as a first exploration and then applies classical random search.
In terms of quasirandom sequence and design of experiments (i.e. other one-shot optimization methods than random search): we use Latin Hypercube Sampling\cite{eglajs,mckay}, Halton\cite{halton60}, Hammersley\cite{hammersley}, and scrambled versions thereof\cite{atanasov,owen,permut1,permut3}; we use the scrambling defined in \cite{permut2}. These algorithms are usually presented in terms of sampling the hypercube $[0,1]^d$; we convert to unbounded domains by applying the inverse cumulative distribution function of the Gaussian for each coordinate. We refer to \cite{qr1} for a more detailed description. For sanity checks and comparison, we include a so-called ``stupid'' optimizer which always outputs $(-1,-1,\dots,-1)$ both as exploration and recommendation.
\item Neural optimization. We use NAS\cite{nas}, adapted for optimizing continuous variables as well\cite{cac}. While this algorithm, internally, uses a gradient, it does not use any gradient of the objective function so that it is a decent competitor for derivative-free optimization.
\end{itemize}

\subsection{Statistical significance}
In all experiments in this section, all x-values correspond to independent runs. Therefore, the frequency at which a curve is below another curve can be used for making statistical inference.
In all experiments, the number between parentheses in the legend is the number of distinct independent experiments used for generating the curve.
In all cases, we display median values.

\subsection{Results on artificial test functions with noise}\label{noise}
We here consider the optimization of functions corrupted by a centered noise, without assuming variance decreasing to zero in the neighborhood of the optimum.
Results are clear: algorithms specifically designed for tackling noise perform best. Noise-managing self-adaptive $(\mu/\mu,\lambda)$-ES\cite{beyernoise}) perform very well; the naive counterpart is less stable; NAS has some merit in moderate dimension; others are weaker (DE variants, DOE, $(1+1)$-ES). The performance of noise-managing ES is particularly visible in Sections \ref{nascooldim5} to \ref{mulambdacool}; it is the only algorithm in the present experiments known for having a simple regret $O(1/budget)$. The performance of the neural net optimization version is visible e.g. in Section \ref{nascooldim5}. Other algorithm sometimes completely fail and diverge far from the optimum, in particular on ill-conditioned functions.
Experiments performed with delay (parallelism) were very similar - the best algorithms are population-based or somehow indifferent to batch-size and are therefore not much impacted by the delay induced by parallelism.

\subsubsection{Dimension 5 - NAS as an interesting competitor to noise-managing-$(\mu/\mu,\lambda)$-ES.}\label{nascooldim5}
First on multimodal functions:

{\center{
\includegraphics[width=.45\textwidth]{tol_figartif2_lognoisy_ackley_function__dim_5__delay_1.png}
\includegraphics[width=.45\textwidth]{tol_figartif2_lognoisy_griewank_function__dim_5__delay_1.png}
\includegraphics[width=.45\textwidth]{tol_figartif2_lognoisy_lunacek_function__dim_5__delay_1.png}
\includegraphics[width=.45\textwidth]{tol_figartif2_lognoisy_rastrigin_function__dim_5__delay_1.png}
\includegraphics[width=.45\textwidth]{tol_figartif2_lognoisy_rosenbrock_function__dim_5__delay_1.png}
\includegraphics[width=.45\textwidth]{tol_figartif2_lognoisy_schwefel_function__dim_5__delay_1.png}
}}

Then on the sphere:

{\center{
\includegraphics[width=.45\textwidth]{tol_figartif2_lognoisy_sphere_function__dim_5__delay_1.png}}}

And then on ill-conditioned unimodal functions:
{\center{
\includegraphics[width=.45\textwidth]{tol_figartif2_lognoisy_alternate_fixedcigar_function__dim_5__delay_1.png}
\includegraphics[width=.45\textwidth]{tol_figartif2_lognoisy_ellipsoid_function__dim_5__delay_1.png}
\includegraphics[width=.45\textwidth]{tol_figartif2_lognoisy_fixedcigar_function__dim_5__delay_1.png}
\includegraphics[width=.45\textwidth]{tol_figartif2_lognoisy_fixeddiscus_function__dim_5__delay_1.png}
}}

\subsubsection{Dimension 20 - noise-managing-$(\mu/\mu,\lambda)$-ES rules.}\label{mulambdacool}
We use large centered translations. Noise-managing self-adaptive $(\mu/\mu,\lambda)$-ES performs best. It consistently outperforms its naive counterpart.

{\center{
\includegraphics[width=.45\textwidth]{tol_figartif2nodoe_lognoisy_ackley_function__dim_20__delay_1.png}
\includegraphics[width=.45\textwidth]{tol_figartif2nodoe_lognoisy_ellipsoid_function__dim_20__delay_1.png}
\includegraphics[width=.45\textwidth]{tol_figartif2nodoe_lognoisy_fixedcigar_function__dim_20__delay_1.png}
\includegraphics[width=.45\textwidth]{tol_figartif2nodoe_lognoisy_fixeddiscus_function__dim_20__delay_1.png}
\includegraphics[width=.45\textwidth]{tol_figartif2nodoe_lognoisy_griewank_function__dim_20__delay_1.png}
\includegraphics[width=.45\textwidth]{tol_figartif2nodoe_lognoisy_lunacek_function__dim_20__delay_1.png}
\includegraphics[width=.45\textwidth]{tol_figartif2nodoe_lognoisy_rastrigin_function__dim_20__delay_1.png}
\includegraphics[width=.45\textwidth]{tol_figartif2nodoe_lognoisy_rosenbrock_function__dim_20__delay_1.png}
}}

\subsection{Results on artificial rotated unimodal functions}\label{manyunimod}
Consistently with the state of the art, for the simple sphere function, the best algorithm (among algorithms robust to compositions with monotonic functions) is the $(1+1)$-ES.

{\center{
\includegraphics[width=.45\textwidth]{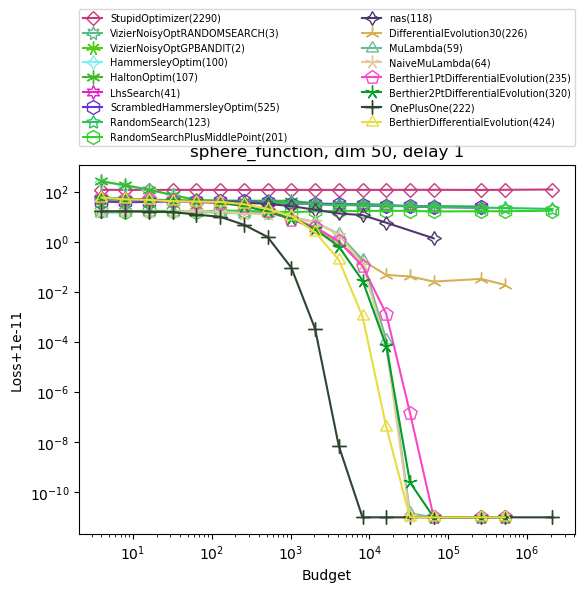}
\includegraphics[width=.45\textwidth]{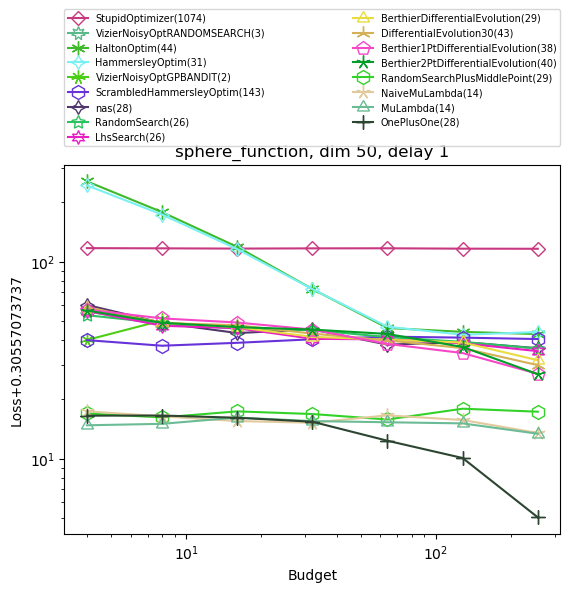}
}}

{\center{\includegraphics[width=.45\textwidth]{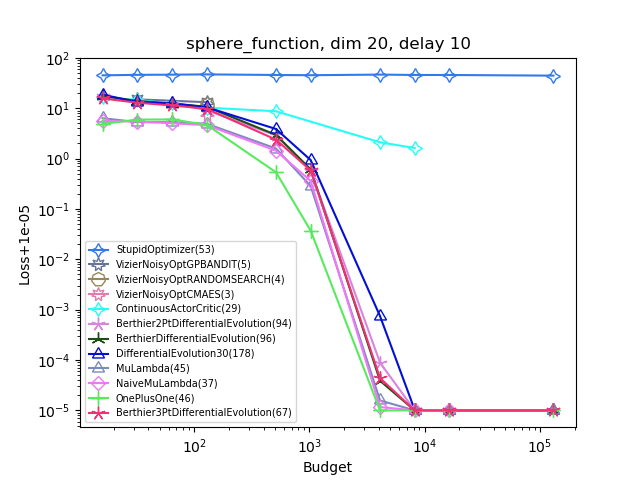}}}

As shown by the third plot, this is still the case when we include a delay (parallelism) of 10.

On partially rotated functions, differential evolution performs quite well:

{\center{
\includegraphics[width=.45\textwidth]{tol_figartif2_logpartially_rot_ackley_function__dim_20__delay_1.png}
\includegraphics[width=.45\textwidth]{tol_figartif2_logpartially_rot_alternate_fixedcigar_function,_dim_20,_delay_1.png}
}}

In some cases, rotational invariance is beneficial; for example when variables correspond to arbitrary axes\cite{rotinv1,rotinv2,rotinv3,rotinv4}. 
However, we might argue that an algorithm which uses the decomposition into axes, by not being invariant by rotation, is more modular than other algorithms\cite{crde}.
The {\em{not}} modular, rotationally invariant DE performs better than other DE when we consider full rotations:

{\em{
\includegraphics[width=.45\textwidth]{tol_figartif2_logrot_alternate_fixedcigar_function__dim_5__delay_1.png}
\includegraphics[width=.45\textwidth]{tol_figartif2_logrot_fixeddiscus_function__dim_5__delay_1.png}
}}

%

\subsection{Results on the non rotated multimodal functions}\label{nonrot}
Without rotation, the original artificial function, though they are not necessarily separable, have some modularity. Unsurprisingly, in this multimodal context, population-based methods (all DE variants, $(\mu/\mu,\lambda)$ variants) perform best. 
Enforcing parallelism does not change much the results - population-based methods already quite parallel by nature.
The $(1+1)$-ES performs well in some cases (in particular for low budget), but sometimes miserably failed. In moderate dimension, continuous actor critic was sometimes competitive.
Experiments here use small translations.

{\center{
\includegraphics[width=.45\textwidth]{tol_figartif2nodoe_logackley_function__dim_20__delay_1.png}
\includegraphics[width=.45\textwidth]{tol_figartif2nodoe_loggriewank_function__dim_20__delay_1.png}
\includegraphics[width=.45\textwidth]{tol_figartif2nodoe_lograstrigin_function__dim_20__delay_1.png}
\includegraphics[width=.45\textwidth]{tol_figartif2nodoe_logackley_function__dim_5__delay_1.png}
\includegraphics[width=.45\textwidth]{tol_figartif2nodoe_loggriewank_function__dim_5__delay_1.png}
\includegraphics[width=.45\textwidth]{tol_figartif2nodoe_lograstrigin_function__dim_5__delay_1.png}
}}

\subsection{Results on artificial rotated multimodal functions, larger translations}\label{rot}
Results vary depending on the budget: the simple $(1+1)$-ES performs well on small budgets, but is outperformed for a larger budget by DE, which is sometimes (not always) outperformed by the $(\mu/\mu,\lambda)$ algorithms, which are sometimes eventually outperformed by NAS.

These figures are for parallelism 10, but results are quite similar in the sequential case.

{\center{
\includegraphics[width=.45\textwidth]{tol_figartif2_logrot_ackley_function__dim_5__delay_10.png}
\includegraphics[width=.45\textwidth]{tol_figartif2_logrot_griewank_function__dim_5__delay_10.png}
\includegraphics[width=.45\textwidth]{tol_figartif2_logrot_lunacek_function__dim_5__delay_10.png}
\includegraphics[width=.45\textwidth]{tol_figartif2_logrot_rastrigin_function__dim_5__delay_10.png}
}}

{\center{
\includegraphics[width=.45\textwidth]{tol_figartif2_logrot_rosenbrock_function__dim_5__delay_10.png}
}}

\subsection{Results on scikit, fast setting}\label{scikitres}
The fast setting implies that the priority is on architectures which are quickly trained.
A short overview is that population-based algorithms are suitable, including $(\mu/\mu,\lambda)$ and DE variants.

\def\toounstable{
\subsubsection{On breast}
The difficult breast dataset leads to poor results with a lot of overfitting. Interestingly, in this setting (retrained architectures, small training and validation set), population-based methods and in particular the quite modular DE perform best. This is consistent with \cite{popbasedflatoptima} - population-based methods prefer stable flat optima.

{\center{
\includegraphics[width=.45\textwidth]{tol_figscikit_logc0_bignn_tuning_for_breast__dim_45__delay_1.png}
\includegraphics[width=.45\textwidth]{tol_figscikit_logc0_nn_tuning_for_breast__dim_9__delay_1.png}
\includegraphics[width=.45\textwidth]{tol_figscikit_logc0_rf_tuning_for_breast__dim_12__delay_1.png}
\includegraphics[width=.45\textwidth]{tol_figscikit_logc0_bignn_tuning_for_breastbig__dim_45__delay_1.png}
\includegraphics[width=.45\textwidth]{tol_figscikit_logc0_nn_tuning_for_breastbig__dim_9__delay_1.png}
\includegraphics[width=.45\textwidth]{tol_figscikit_logc0_rf_tuning_for_breastbig__dim_12__delay_1.png}
}}
}
\subsubsection{With the small datasets.}\label{small}
With small datasets, we often see increasing test loss while we optimize the validation loss. This is overfitting.
We get a good performance for all variants of DE,
with less overfitting. No other algorithm could compete in terms of stability. Our interpretation is that the combination in DE smoothly interpolate fitness values and generate stable individuals. 
This is consistent with \cite{popbasedflatoptima}; but tells us a bit more - DE here outperformed clearly other population-based methods, even those that perform very well in noisy setting or difficult rotated/multimodal testbeds (noise-managing ES, NAS).
In terms of high level chosen architectures, RF is competitive except for digits (images).

{\center{
\includegraphics[width=.45\textwidth]{tol_figscikit_logc0_bignn_tuning_for_diabetes__dim_45__delay_1.png}
\includegraphics[width=.45\textwidth]{tol_figscikit_logc0_bignn_tuning_for_iris__dim_45__delay_1.png}
\includegraphics[width=.45\textwidth]{tol_figscikit_logc0_nn_tuning_for_diabetes__dim_9__delay_1.png}
\includegraphics[width=.45\textwidth]{tol_figscikit_logc0_rf_tuning_for_diabetes__dim_12__delay_1.png}
}}

\subsubsection{With the full datasets.}
We get a good performance
for all variants of DE. DE, which performed well on separable or partially rotated problems but not that much (except the rotationally invariant version) on rotated problems, is effective on AS. TODO more on this: results in the non-rotated case for the rotinv DE ?

{\center{
\includegraphics[width=.45\textwidth]{tol_figscikit_logc0_bignn_tuning_for_digitsbig__dim_45__delay_1.png}
\includegraphics[width=.45\textwidth]{tol_figscikit_logc0_nn_tuning_for_diabetesbig__dim_9__delay_1.png}}}

{{\center
\includegraphics[width=.45\textwidth]{tol_figscikit_logc0_nn_tuning_for_digitsbig__dim_9__delay_1.png}
}}
}

\section{Conclusion}
We have shown mathematically that the test-based population-size adaptation from \cite{beyernoise} can be applied in other contexts, namely multimodal optimization, in which it  escapes plateaus.
It does not escape convex local minima, but experimentally finite machine precision is enough for naturally ensuring some kind of restart.
\otc{The method can be applied to other optimization algorithms.}

\subsection*{Further work}
We got good results in parallel settings, with moderate numbers of iterations. Results are less satisfactory for cases with large budget. We guess that applying the same TBPSA mechanism on top of other algorithms (CMA, CMSA, DE, PSO) might be beneficial.
A limited precision in the test might also facilitate the detection of convex local minima.

\bibliographystyle{ACM-Reference-Format}
\bibliography{lsca,doe}

\end{document}